\documentclass[11pt,reqno]{amsart}
\usepackage{amsmath,amsthm,amsfonts,amssymb,xcolor,comment,hyperref,comment}
\usepackage{graphicx}
\usepackage{setspace}
\usepackage{mathrsfs}
\usepackage{hyperref,bbm}
\usepackage[letterpaper,margin=1.3in]{geometry}

\usepackage[shortlabels]{enumitem}

\usepackage{comment}

\usepackage[
    backend=biber,
    style=numeric,
    doi=false,isbn=false,url=false,eprint=false
    ]{biblatex}
\renewcommand{\epsilon}{\varepsilon}

\renewcommand{\tilde}{\widetilde}

\def\XXint#1#2#3{{\setbox0=\hbox{$#1{#2#3}{\int}$}
     \vcenter{\hbox{$#2#3$}}\kern-.5\wd0}}

\newtheorem{theorem}{Theorem}[section]
\newtheorem*{theorem*}{Theorem}
\newtheorem{proposition}[theorem]{Proposition}
\newtheorem{lemma}[theorem]{Lemma}
\newtheorem{corollary}[theorem]{Corollary}

\newtheorem{example}[theorem]{Example}

\theoremstyle{definition}
\newtheorem{definition}[theorem]{Definition}

\theoremstyle{remark}
\newtheorem{remark}[theorem]{Remark}

\numberwithin{equation}{section}

\newcommand{\bbG}{\mathbb{G}}
\newcommand{\bbR}{\mathbb{R}}

\title[Quotient maps, pliability, and Lusin approximation]{Quotient maps, pliability, and Lusin approximation in Carnot groups}

\author[Gareth Speight]{Gareth Speight}
\address{Department of Mathematical Sciences, University of Cincinnati, 2815 Commons Way, Cincinnati, OH 45221, United States}
\email{Gareth.Speight@uc.edu}

\author[Scott Zimmerman]{Scott Zimmerman}
\address{Department of Mathematics, 
The Ohio State University,
100 Math Tower, 231 West 18th Avenue, Columbus, OH 43210, United States}
\email{zimmerman.416@osu.edu}


\allowdisplaybreaks

\begin{document}

\begin{abstract}
In Carnot groups, we show that pliability and the Lusin approximation property are preserved by quotient maps, and we deduce several applications. For instance, we show that, in free Carnot groups of step at least 3 with at least $2$ generators, every horizontal direction is non-pliable, the Lusin property fails in a strong sense, and there is a Lipschitz curve that is purely unrectifiable with respect to smooth horizontal curves. We also give a previously unpublished example of a pliable Carnot group, namely, the group of upper triangular, unipotent matrices.
\end{abstract}


\date{\today}

\maketitle

\section{Introduction}

Lusin's classical theorem asserts that any measurable function on $\mathbb{R}$ coincides with a continuous function outside a set of arbitrarily small measure. More recently,
similar Lusin approximation results have played 
an important role in the fields of geometric measure theory and geometric function theory. 
If one imposes stronger hypotheses on the starting map, one can approximate by maps with higher regularity. For example, every absolutely continuous map can be approximated by $C^{1}$ mappings.
A standard strategy for constructing these Lusin-type approximations relies on the 
Whitney Extension Theorem \cite{Whi34}.
Given a map of low regularity, one identifies a large subset of the domain 
on which the hypotheses of Whitney's theorem holds, and this result is then applied to create a smooth map that differs from the original one on a small set.

In recent years, these questions have been 
studied in non-Euclidean settings, in particular in Carnot groups and other sub-Riemannian manifolds \cite{BLU07, Pan89}. A Carnot group $\mathbb{G}$ is a connected, simply connected Lie group whose Lie algebra $\mathfrak{g} = \mathfrak{g}_1 \oplus \dots \oplus \mathfrak{g}_s$ is stratified. The first layer $\mathfrak{g}_1$ defines a horizontal distribution, and curves whose velocity vectors lie in $\mathfrak{g}_1$ almost everywhere are called \emph{horizontal curves}. Such curves represent the directions through which one may travel in such a manifold. 
If these velocity vectors vary continuously in $\mathfrak{g}_1$, we say that the curve is of class $C^1_H$.
Because the geometry of a sub-Riemannian manifold is defined by its horizontal directions, classical Euclidean extension and approximation techniques often fail or require delicate sub-Riemannian adjustments.

Juillet and Sigalotti introduced the notion of \emph{pliability} for vectors in the horizontal layer of a Carnot group \cite{JS17}.
A vector $V \in \mathfrak{g}_1$ is pliable if arbitrarily small 
$C^1_H$ perturbations of the curve $t \mapsto \exp(tV)$ can 
always reach an open neighborhood of the endpoint $\exp(V)$ in $\mathbb{G}$. Juillet and Sigalotti proved 
that a Whitney Extension Theorem for $C^1_H$ curves holds in a Carnot group $\mathbb{G}$ if and only if every vector in the horizontal layer of $\mathbb{G}$ is pliable.  These properties then imply a $C^{1}_{H}$ Lusin approximation for Lipschitz curves \cite{JS17, LDS16, Spe16, Zim18}. 

In step 2 Carnot groups, all horizontal vectors are pliable. However, non-pliable vectors emerge in Carnot groups of step at least 3. For example, in the step 3 Engel group $\mathbb{E}$, the horizontal vector $E_2$ is non-pliable, whereas all other non-zero horizontal directions are pliable \cite{JS17, SZ25}. In addition, it can be shown that the Engel group does not admit $C^{1}_{H}$ Lusin approximation for Lipschitz curves \cite{Spe16}. However, it does admit a weak Lusin approximation property, where every Lipschitz curve meets a $C^{1}$ horizontal curve in a set of strictly positive measure. This was shown in  \cite{SZ25}, as an application of the fact that pliability in a subset of horizontal directions is still sufficient to guarantee the existence of $C^1_H$ Whitney extensions and Lusin approximations for curves in those directions. The authors showed in \cite{SZ26} that the free Carnot group $\mathbb{F}_{2,3}$ of step 3 with 2 generators exhibits even more pathological behavior than the Engel group. Indeed, there exists a Lipschitz curve whose intersection with  {\em any} $C^1_H$ curve has 1-dimensional Hausdorff measure 0. That is, Lipschitz curves in $\mathbb{F}_{2,3}$ may not necessarily be approximated by smooth curves.

In the present paper we investigate two natural questions. Firstly, we study how pliability and Lusin approximation are preserved by suitable mappings by proving that these properties are preserved by quotient maps between the Lie algebras $\mathfrak{g}$ and $\mathfrak{h}$ associated to Carnot groups. A quotient map (Definition \ref{quotientmap}) is a Lie algebra homomorphism that maps the horizontal layer $\mathfrak{g}_1$ of $\mathfrak{g}$ 
onto the horizontal layer $\mathfrak{h}_1$ of $\mathfrak{h}$. 
In summary, we prove the following result.

\begin{theorem}
\label{t-combine}
Suppose $\mathbb{G}, \mathbb{H}$ are Carnot groups, 
$\mathfrak{g}, \mathfrak{h}$ are their Lie algebras with horizontal layers $\mathfrak{g}_1,\mathfrak{h}_1$ respectively, and
$\Phi \colon \mathfrak{g} \to \mathfrak{h}$ is a quotient map.
\begin{enumerate}
\item If $\mathfrak{B}\subset \mathfrak{h}_{1}$ and $\mathbb{G}$ has the $C^1_H$ Lusin property on $\Phi^{-1}(\mathfrak{B}) \subset \mathfrak{g}_1$ , then $\mathbb{H}$ has the $C^1_H$ Lusin property on $\mathfrak{B}$.
\item If $\mathbb{G}$ has the weak $C^1_H$ Lusin property, then $\mathbb{H}$ has the weak $C^1_H$ Lusin property.
\item \label{enum:pliable} If $X\in \mathfrak{g}_1$ is pliable, then $\Phi(X)\in \mathfrak{h}_1$ is pliable.
\end{enumerate}
\end{theorem}

We note that property (\ref{enum:pliable}) in Theorem~\ref{t-combine}
was mentioned in a footnote on page 1644 of \cite{JS17} as an observation attributed to E. Le Donne. In this paper, we provide the precise details of the argument justifying this statement.

A second natural question is to what extent the pathological curves constructed in $\mathbb{F}_{2,3}$ exist in other Carnot groups. We show that these results generalize to many more Carnot groups, namely all free Carnot groups $\mathbb{F}_{r,s}$ of step $s \geq 3$ with $r \geq 2$ generators. 
Write $\mathfrak{f}_{r,s}$ to denote the Lie algebra of $\mathbb{F}_{r,s}$.

\begin{theorem}
    \label{t-purelyunrec}
        Fix integers $r \geq 2$ and $s \geq 3$.
        \begin{enumerate}
            \item The only pliable vector in $\mathfrak{f}_{r,s}$ is $0$.
            \item There is a Lipschitz curve $\gamma\colon[0,1] \to \mathbb{F}_{r,s}$ such that
    $\mathcal{H}^1(\operatorname{Im}(\gamma) \cap \operatorname{Im}(\Gamma)) = 0$ for any $C^1_H$ curve $\Gamma$ in $\mathbb{F}_{r,s}$.
        \end{enumerate}
\end{theorem}

To prove Theorem \ref{t-purelyunrec} we show that there is a quotient map from $\mathbb{F}_{r,s}$ to $\mathbb{F}_{2,3}$ whenever $r\geq 2, s\geq 3$. It was shown in \cite{SZ26} that $\mathbb{F}_{2,3}$ does not have the weak $C^{1}_H$ Lusin property. Hence by Theorem \ref{t-combine}, we will be able to conclude that $\mathbb{F}_{r,s}$ does not have the weak $C^{1}_{H}$ Lusin property. Combining this with \cite[Theorem 5.1]{SZ26} then yields Theorem \ref{t-purelyunrec}. Again using Theorem \ref{t-combine} and similar techniques, we also study pliable vectors in model Filiform groups (Corollary \ref{c-fil}) and unipotent matrix groups (Corollary \ref{c-nipotentexample}).
In particular, we show that the group of $(s + 1) \times (s+1)$ unipotent matrices is another example of a pliable Carnot group of step $s$ for any $s$.

The paper is organized as follows. In Section~\ref{s-carnot}, we recall preliminaries on Carnot groups, 
horizontal curves, pliability, and Lusin properties. In Section~\ref{s-directional}, we establish that quotient maps preserve Lusin approximation properties. In Section~\ref{s-pliablepush}, we show that quotient maps push forward pliable vectors. This proves Theorem \ref{t-combine}. Finally, in Section~\ref{s-applications}, we construct several quotient maps and derive some applications of Theorem \ref{t-combine}. In particular, we prove Theorem \ref{t-purelyunrec}

\medskip

\textbf{Acknowledgements:} G. Speight was supported by the National Science Foundation under Award No. 2348715.

\section{Preliminaries}
\label{s-carnot}

\subsection{Carnot Groups}

We briefly recall basic terminology involving Carnot groups. A Lie algebra is a vector space $\mathfrak{g}$ equipped with a Lie bracket $[\cdot,\cdot] \colon \mathfrak{g}\times \mathfrak{g}\to \mathfrak{g}$ that is bilinear, satisfies $[X,X]=0$, and satisfies the Jacobi identity $[X,[Y,Z]]+[Y,[Z,X]]+[Z,[X,Y]]=0$. The Lie algebra associated to a Lie group is the space of left invariant vector fields equipped with the Lie bracket $[X,Y]$ defined on smooth functions as $[X,Y](f)=X(Y(f))-Y(X(f))$.

\begin{definition}\label{Carnot}
A Lie algebra $\mathfrak{g}$ is {\em nilpotent} and {\em stratified} if
it admits a decomposition as a direct sum of subspaces of the form
\[\mathfrak{g}=\mathfrak{g}_{1}\oplus \mathfrak{g}_{2}\oplus \cdots \oplus \mathfrak{g}_{s}\]
such that $\mathfrak{g}_{i}=[\mathfrak{g}_{1},\mathfrak{g}_{i-1}]$ for any $i=2, \ldots, s$, and $[\mathfrak{g}_{1},\mathfrak{g}_{s}]=0$. The subspace $\mathfrak{g}_{1}$ is called the \emph{horizontal} or {\em first layer}. 

A Lie group $\bbG$ is a {\em Carnot group} if it is connected, simply connected, and its Lie algebra $\mathfrak{g}$ is nilpotent and stratified.
\end{definition}

Let $\bbG$ be a Carnot group. We fix a basis $\{X_1,\dots,X_r\}$ of its horizontal layer. This equips $\mathfrak{g}_1$ with an inner product $\langle\cdot,\cdot\rangle_{\mathfrak{g}_1}$ for which $\{X_1,\dots,X_r\}$ is orthonormal. We then define the length $\Vert \cdot \Vert_{\mathfrak{g}_1}$ of a vector $V = \sum_{i=1}^r v_i X_i \in \mathfrak{g}_1$ by $\Vert V \Vert_{\mathfrak{g}_1} = \sqrt{\langle V,V \rangle_{\mathfrak{g}_1}} = \left( \sum_{i=1}^r v_i^2\right)^{1/2}$.

\begin{remark}
In any Carnot group $\mathbb{G}$ we adopt the following notation. Let $m_i:=\dim(\mathfrak{g}_i)$, $h_i:=m_1+\dots +m_i$ for $1\leq i\leq s$, and $h_{0}:=0$. A basis $X_1,\dots, X_n$ of $\mathfrak{g}$ is \emph{adapted to the stratification} if $X_{h_{i-1}+1},\dots, X_{h_{i}}$ is a basis of $\mathfrak{g}_i$ for $1\leq i \leq s$. We fix such a basis. Let $r$ be the dimension of the horizontal layer i.e. $r:=m_{1}$.
\end{remark}

\begin{definition}
The map $\exp \colon \mathfrak{g}\to \mathbb{G}$ is defined by $\exp(X)=\gamma(1)$ where $\gamma\colon \bbR \to \bbG$ is the unique one-parameter subgroup of 
$\bbG$ whose tangent vector at the identity is equal to $X$.
This \emph{exponential map} is a diffeomorphism between $\mathfrak{g}$ and $\mathbb{G}$. 
\end{definition}

We can identify $\mathbb{G}$ with $\mathbb{R}^{n}$ using the following diffeomorphism:
\[ (x_{1}, \ldots, x_{n})\in \mathbb{R}^{n} \longleftrightarrow \exp(x_{1}X_{1}+\cdots +x_{n}X_{n})\in \mathbb{G}.\]
This identification is known as \emph{exponential coordinates of the first kind}.
It induces a Lie group operation on $\bbR^{n}$ according to the Baker-Campbell-Hausdorff formula  \cite[Definition~2.2.11]{BLU07}.

A curve $\gamma \colon [a,b]\to \bbG$ is \emph{absolutely continuous} if it is absolutely continuous as a curve into $\mathbb{R}^{n}$. This implies  the velocity $\dot{\gamma}(t)$ exists for almost every $t\in [a,b]$.
Since the choice of coordinates equips $\bbR^n$ with a Lie group structure with Lie algebra $\mathfrak{g}$, we may view any $X \in \mathfrak{g}$ as a left invariant vector field on $\mathbb{R}^n$. That is, we can write 
$\dot{\gamma}(t)=\sum_{j=1}^{n}u_{j}(t)X_{j}(\gamma(t))$ for some controls $u_{1}, \ldots, u_{n}\in L^{1}[a,b]$.

\begin{definition}
An absolutely continuous curve $\gamma\colon [a,b]\to \bbG$ is \emph{horizontal} if we can write 
\begin{equation}\label{dotgamma}
\dot{\gamma}(t)=\sum_{j=1}^{r}u_{j}(t)X_{j}(\gamma(t)) \qquad \mbox{almost everywhere},
\end{equation}
i.e. if $u_{r+1} = \cdots = u_n = 0$ a.e. We define the \emph{length} of such a curve by $\ell_{\bbG}(\gamma):=\int_{a}^{b} |u|$.
\end{definition}

For any $t$ at which $\dot{\gamma}(t)$ exists, 
we will write 
\begin{equation}
\label{e-derivativeversions}
\gamma'(t) := 
(dL_{\gamma(t)^{-1}})_{\gamma(t)} \dot{\gamma}(t)
=
\sum_{j=1}^{n}u_{j}(t)X_{j} \in \mathfrak{g}_1
\end{equation}
where $L_{\gamma(t)^{-1}}\colon \mathbb{G} \to \mathbb{G}$ is left multiplication by $\gamma(t)^{-1}$, i.e. $L_{\gamma(t)^{-1}}(p) = \gamma(t)^{-1}*p$, and $(dL_{\gamma(t)^{-1}})_{\gamma(t)}$ is its Jacobian matrix at $\gamma(t)$.
In doing so we view $\gamma'$ as a map into $\mathfrak{g}$, which we identify with the tangent space to $\mathbb{G}$ at the identity.
This gives the following alternate definition of horizontality.

\begin{definition}
An absolutely continuous curve $\gamma\colon [a,b]\to \bbG$ is \emph{horizontal} if $\gamma' \in \mathfrak{g}_1$ almost everywhere.
\end{definition}

\begin{definition}
The \emph{Carnot-Carath\'{e}odory distance (CC distance)} between any two points $x, y\in \bbG$ is defined by
\[d(x,y):=\inf \{\ell_{\bbG}(\gamma)\colon \gamma \mbox{ is a horizontal curve joining } x \mbox{ and }y\}.\]
\end{definition}

Note that the CC distance is left-invariant 
i.e. $d(r*p,r*q) = d(p,q)$
for any $p,q,r \in \bbG$.
The CC distance induces on $\mathbb{G}$ the same topology as the Euclidean distance. 
However, it is not bi-Lipschitz equivalent to the Euclidean distance.

\subsection{Pliability}

Suppose $\mathbb{G}$ is a Carnot group with Lie algebra $\mathfrak{g}$.

\begin{definition}
    Given an interval $I \subset \bbR$, a curve $\gamma \colon I \to \bbG$ is $C^1_H$ if it is $C^1$ (i.e. continuously differentiable) as a curve in $\bbR^n$ and if it is a horizontal curve. Equivalently, $\gamma$ is $C^1_H$ if $\gamma' \colon I \to \mathfrak{g}_1$ exists and is continuous.
\end{definition}

We denote by $C_{H}^{1}(I,\bbG)$ the set of all $C^1_H$ curves $\gamma \colon I \to \bbG$. 
If $I=[a,b]$, the space $C_{H}^{1}([a,b],\bbG)$ can be equipped with a natural metric:
\[ d_{C^1_H}(\gamma_{1},\gamma_{2}):=\max \left( \sup_{t\in [a,b]} d_{\mathbb{G}}(\gamma_{1}(t),\gamma_{2}(t)), \ \sup_{t\in [a,b]} \|\gamma_{2}'(t)-\gamma_{1}'(t)\|_{\mathfrak{g}_1} \right).
\]
We will sometimes use the notation $\|\gamma_{2}'-\gamma_{1}'\|_{\infty}$ to denote $\sup_{t\in [a,b]} \|\gamma_{2}'(t)-\gamma_{1}'(t)\|_{\mathfrak{g}_1}$. 
If working in coordinates, we will occasionally use $|\cdot|$ to denote the Euclidean norm on $\mathbb{G}$. 

The following geometric property of curves was introduced in \cite{JS17}.
The following definition was not the original one, but their equivalence was verified in \cite[Proposition 3.7]{JS17}.

\begin{definition}\label{pliability}
A vector $V\in \mathfrak{g}_{1}$ is pliable if, for every neighborhood $\mathcal{V}$ of the curve $ t\mapsto \exp(tV)$ in the space $C_{H}^{1}([0,1],\bbG)$, the set
\[ \{\beta(1): \beta \in \mathcal{V}, (\beta, \beta')(0)=(0_{\mathbb{G}}, V) \} \]
is a neighborhood of $\exp(V)$.
\end{definition}

We say that $\mathbb{G}$ is {\em pliable} if every $V \in \mathfrak{g}_1$ is pliable.
The importance of this property is its equivalence to the Whitney Extension Property for $C^1_H$ curves in Carnot groups, which was the main result of \cite{JS17}. It is also proven in \cite{JS17} that all step two Carnot groups are pliable and that there exist pliable Carnot groups of arbitrarily high step.

\subsection{Lusin approximation}

Suppose $\mathbb{G}$ is a Carnot group with Lie algebra $\mathfrak{g}$.

\begin{definition}\label{defLusin}
Let $\mathfrak{A}\subset \mathfrak{g}_1$. We say that $\mathbb{G}$ has the {\em $C^1_H$ Lusin property on $\mathfrak{A}$} if, for any $\varepsilon>0$ and any horizontal curve $\gamma \colon [0,1] \to \mathbb{G}$, there is a $C^1_H$ curve $\Gamma \colon [0,1]\to \mathbb{G}$ such that 
\[
m\{ t \in [0,1] : \gamma'(t) \in \mathfrak{A} 
\text{ and }
\Gamma(t) \neq \gamma(t)\} 
<
\varepsilon.
\]
If $\mathfrak{A} = \mathfrak{g}_{1}$, this inequality is simply
\[
m\{t \in [0,1] : \Gamma(t) \neq \gamma(t) \} < \varepsilon.
\]
Hence we recover the usual $C^1_H$ Lusin property.
\end{definition}

Intuitively, Definition \ref{defLusin} means that the parts of horizontal curves which travel in direction of $\mathfrak{A}$ can be well approximated by $C^{1}$ horizontal curves.

\begin{definition}\label{defweakLusin}
We say that $\mathbb{G}$ has the {\em weak $C^1_H$ Lusin property} if, for any horizontal curve $\gamma\colon [0,1] \to \mathbb{G}$, there is a $C^1_H$ curve $\Gamma\colon [0,1]\to \mathbb{G}$ such that 
\[
m\{ t \in [0,1] : 
\Gamma(t) \neq \gamma(t)\} 
<
1.
\]
\end{definition}

In other words, Definition \ref{defweakLusin} means that every horizontal curve must agree with a $C^1_H$ curve on a set of positive measure, but we have no quantitative control over the measure of this set.

It follows from \cite{JS17} that all pliable Carnot groups have the $C^1_H$ Lusin property, while the Engel group does not \cite{Spe16}. The Engel group does, however, have the weak $C^1_H$ Lusin property \cite{SZ25}. The free Carnot group of rank 2 and step 3 does not even have the weak $C^1_H$ Lusin property \cite{SZ26}.

\subsection{Lie algebra homomorphisms}

The following definitions are from \cite[Section 2.1.8]{BLU07}.

\begin{definition}
Let $(\bbG, \ast_{1})$, $(\mathbb{H}, \ast_{2})$ be Lie groups. A map $\varphi\colon \bbG\to \mathbb{H}$ is a \emph{Lie group homomorphism} if it is $C^{\infty}$ and
\[\varphi(x\ast_{1} y)=\varphi(x)\ast_{2} \varphi(y) \qquad \mbox{for all }x,y\in \bbG.\]
A map is a {\em Lie group isomorphism} if it is a Lie group homomorphism and a diffeomorphism of differentiable manifolds.
\end{definition}

\begin{definition}
Let $(\mathfrak{g},[\cdot,\cdot]_{1})$, $(\mathfrak{h},[\cdot,\cdot]_{2})$ be Lie algebras. A map $\Phi\colon \mathfrak{g}\to \mathfrak{h}$ is a \emph{Lie algebra homomorphism} if it is linear and
\[\Phi([X,Y]_{1})=[\Phi(X),\Phi(Y)]_{2} \qquad \mbox{ for all }X, Y\in \mathfrak{g}.\]
A map is a {\em Lie algebra isomorphism} if it is a bijective Lie algebra homomorphism.
\end{definition}

The following result is fundamental in Lie theory.

\begin{theorem}[Lee, 2013, Thm 20.19]
\label{t-fundamental}
    Suppose $\mathbb{G},\mathbb{H}$ are Carnot groups with Lie algebras $\mathfrak{g},\mathfrak{h}$ respectively.
    Given a Lie algebra homomorphism $\Phi\colon \mathfrak{g}\to \mathfrak{h}$, 
    the map $F\colon \mathbb{G} \to \mathbb{H}$ 
    defined as
    \[
    F (\exp_{\mathbb{G}}(X)) = \exp_{\mathbb{H}}(\Phi(X))
    \]
    is a 
    Lie group homomorphism, and $dF_{e_\mathbb{G}} = \Phi$.
\end{theorem}

We consider the following class of maps.

\begin{definition}\label{quotientmap}
    A Lie algebra homomorphism $\Phi \colon \mathfrak{g} \to \mathfrak{h}$ between stratified, nilpotent Lie algebras is called a {\em quotient map} if $\Phi(\mathfrak{g}_1) = \mathfrak{h}_1$.
\end{definition}

The following is simply a rewording of the standard chain rule for differentials 
\cite[Proposition~3.5]{Lee13}.
It will follow from our definition \eqref{e-derivativeversions} of $\gamma'$.

\begin{proposition}
\label{p-chainrule}
Suppose $\Phi\colon \mathfrak{g} \to \mathfrak{h}$ is a Lie algebra homomorphism between stratified, nilpotent Lie algebras and $F\colon \mathbb{G}\to\mathbb{H}$ is its induced Lie group homomorphism.
If $\gamma \colon [0,1] \to \mathbb{G}$ is a horizontal curve, then, for almost every $t \in [0,1]$,
\begin{align*}
    (F \circ \gamma)'(t) = \Phi(\gamma'(t)).
\end{align*}
\end{proposition}

\begin{proof}
Note that $L_{F(p)} \circ F = F \circ L_p$ for any $p \in \mathbb{G}$ since $F$ is a Lie group homomorphism. 
Fix $t \in [0,1]$ at which $\gamma$ is differentiable.
Thus, writing $p = \gamma(t)$,
\eqref{e-derivativeversions} gives
\begin{align*}
    (F \circ \gamma)'(t) 
    = 
    (dL_{F(p)^{-1}})_{F(p)} \tfrac{d}{dt} (F \circ \gamma)(t)
    &=(dL_{F(p)^{-1}})_{F(p)}
    dF_{p}(\dot{\gamma}(t))
    \\
    &= d(L_{F(p)^{-1}} \circ F)_{p}
    (\dot{\gamma}(t))
    \\
    &= d(F \circ L_{p^{-1}})_{p}
    (\dot{\gamma}(t))
    \\
    &=
    dF_{e_\mathbb{G}}(dL_{p^{-1}})_{p}(\dot{\gamma}(t))
    = \Phi(\gamma'(t))
    \end{align*}
    where the third and fifth equalities here follow from the standard chain rule for differentials \cite[Proposition~3.5]{Lee13}.
\end{proof}

\begin{corollary}
    Suppose $\Phi\colon \mathfrak{g} \to \mathfrak{h}$ is a quotient map and $F\colon \mathbb{G} \to \mathbb{H}$ is its induced Lie group homomorphism. 
    If $\alpha\colon [0,1] \to \mathbb{G}$ is horizontal, then $F \circ \alpha$ is horizontal in $\mathbb{H}$.
    If $\alpha \in C_{H}^{1}([0,1],\mathbb{G})$, then $F \circ \alpha \in C^1_H([0,1],\mathbb{H})$.
\end{corollary}

\section{Quotient maps preserve Lusin approximations}
\label{s-directional}

In this section we observe that Lusin properties of Lipschitz curves in Carnot groups are preserved under quotient maps. 
First, we record the following lemma proving that horizontal curves can be ``lifted'' by quotient maps so that the length of the lift is controlled by the length of the original curve.
A similar result was proven in Theorem~4.1 of \cite{LDS16}.

\begin{lemma}
    \label{e-PicardLindeloff}
    Fix Carnot groups $\mathbb{G}, \mathbb{H}$ and
    a quotient map $\Phi\colon  \mathfrak{g} \to \mathfrak{h}$. We denote by $F\colon \mathbb{G} \to \mathbb{H}$ the Lie group homomorphism induced by 
    $\Phi$, i.e., $dF_{e_\mathbb{G}} = \Phi$.
    There is a constant $C>0$ such that,
    for any horizontal curve $\gamma\colon [0,1] \to \mathbb{H}$ with $\gamma(0) = 0_{\mathbb{H}}$, there is a horizontal curve $\alpha\colon [0,1] \to \mathbb{G}$ satisfying $F \circ \alpha = \gamma$ everywhere, $\Phi(\alpha') = \gamma'$ almost everywhere, 
    and
    \begin{equation}
        \label{e-lengthboundlemma}
        L_\mathbb{G}(\alpha) \leq C L_\mathbb{H}(\gamma).
    \end{equation}
    \end{lemma}

\begin{proof}
    Since $\Phi$ is a quotient map, the restriction $\Phi|_{(\ker \Phi)^\perp \cap \mathfrak{g}_1}$ is an isomorphism of $(\ker \Phi)^\perp \cap \mathfrak{g}_1$ onto 
    $\mathfrak{h}_1$.
    Write 
    \[
    A := (\Phi|_{(\ker \Phi)^\perp \cap \mathfrak{g}_1} )^{-1} \colon 
    \mathfrak{h}_1
    \to (\ker \Phi)^\perp \cap \mathfrak{g}_1.
    \]
    Suppose $\gamma\colon [0,1] \to \mathbb{H}$ is a horizontal curve with $\gamma(0) = 0_\mathbb{H}$.
    The IVP
    \[
    \begin{cases}
        \alpha'(t) = A(\gamma'(t)) \mbox{ for a.e. }t\in [0,1]\\
        \alpha(0)=0_\mathbb{G}
    \end{cases}
    \]
    has a unique solution $\alpha\colon [0,1] \to \mathbb{G}$
    (see, for example, \cite[Theorem~A.3]{rifford}).
    Moreover,
    $\alpha$ is horizontal since $\alpha' = A(\gamma') \in (\ker \Phi)^\perp \cap \mathfrak{g}_1$.
    In particular, $\Phi\circ (\alpha') = \Phi \circ A \circ (\gamma') = \gamma'$, and by Proposition \ref{p-chainrule},
    \[
    (F \circ \alpha)'(t) = \Phi(\alpha'(t)) 
    = \gamma'(t)
    \]
    for a.e. $t \in [0,1]$. Since $F(\alpha(0)) = F(0_\mathbb{G}) = 0_\mathbb{H}$, it follows that
    $F \circ \alpha = \gamma$ a.e.
    Finally,
    \[
    L_\mathbb{G}(\alpha) = \int_0^1 \|\alpha'\|_{\mathfrak{g}_1} = \int_0^1 \|A \circ (\gamma')\|_{\mathfrak{g}_1} 
    \leq \Vert A \Vert \int_0^1\|\gamma'\|_{\mathfrak{h}_1}  = \Vert A \Vert L_\mathbb{H}(\gamma)
    \]
    where $\Vert A \Vert$ is the operator norm of the linear map $A$.
\end{proof}

While any Carnot group that has the $C^1_H$ Lusin property (such as every step 2 group) must have the weak $C^1_H$ Lusin property, Theorem~1.2 from \cite{SZ26} shows that there are Carnot groups that do not.
We now show that both properties are preserved by quotient maps.

\begin{proposition}
\label{t-surjective}
Suppose $\mathbb{G}, \mathbb{H}$ are Carnot groups and $\Phi\colon \mathfrak{g} \to \mathfrak{h}$ is a quotient map.
\begin{enumerate}
\item If $\mathbb{G}$ has the weak $C^1_H$ Lusin property, then $\mathbb{H}$ has the weak $C^1_H$ Lusin property.
\item If $\mathfrak{B}\subset \mathfrak{h}_{1}$ and $\mathbb{G}$ has the $C^1_H$ Lusin property on 
$\Phi^{-1}(\mathfrak{B})  \cap \mathfrak{g}_1$,
then $\mathbb{H}$ has the $C^1_H$ Lusin property on $\mathfrak{B}$. 
\end{enumerate}
\end{proposition}

\begin{proof}
Let $F:\mathbb{G} \to \mathbb{H}$ be the Lie group homomorphism 
associated with $\Phi$ that is
guaranteed by Theorem \ref{t-fundamental}. Assume that $\gamma\colon [0,1] \to \mathbb{H}$ is horizontal, and assume without loss of generality that $\gamma(0)=0$.
By Lemma~\ref{e-PicardLindeloff}, there is a horizontal curve $\alpha\colon [0,1] \to \mathbb{G}$ such that
$\gamma = F \circ \alpha$
and $\Phi(\alpha'(t)) = \gamma'(t)$ for a.e. $t \in [0,1]$.

Suppose $\beta\colon [0,1] \to \mathbb{G}$ is any $C^1_H$ curve.
Then the curve $\Gamma = F \circ \beta\colon [0,1] \to \mathbb{H}$ is $C^1_H$, and, 
for all $t \in [0,1]$ such that $\alpha(t) = \beta(t)$, we have
\[
\Gamma(t) = F(\beta(t)) = F(\alpha(t)) = \gamma(t).
\]
Therefore,
$\{ t \in [0,1] : 
\Gamma(t) \neq \gamma(t)\}
\subset
\{ t \in [0,1] : 
\beta(t) \neq \alpha(t)\}
$.

With this observation we now establish (1) and (2). Suppose $\mathbb{G}$ has the weak $C^1_H$ Lusin property. 
Then, in the argument above, we may additionally assume that $\beta$ satisfies
\[ m\{t\in [0,1]:\beta(t)\neq \alpha(t)\}<1,\]
from which it follows that $m\{t\in [0,1]:\Gamma(t)\neq \gamma(t)\}<1$.
This proves (1).

Next suppose $\mathbb{G}$ has the $C^1_H$ Lusin property on $\Phi^{-1}(\mathfrak{B})  \cap \mathfrak{g}_1$. Given $\varepsilon>0$, we may assume that $\beta$ satisfies
\[
m\{ t \in [0,1] : \alpha'(t) \in \Phi^{-1}(\mathfrak{B})  \cap \mathfrak{g}_1 
\text{ and }
\beta(t) \neq \alpha(t)\} 
<
\varepsilon.
\]
If $\gamma'(t) \in \mathfrak{B}$ and $\Phi(\alpha'(t)) = \gamma'(t)$, we have $\alpha'(t)\in \Phi^{-1}(\mathfrak{B})$. Since $\Phi(\alpha'(t)) = \gamma'(t)$ and $\alpha'(t)\in \mathfrak{g}_{1}$ for a.e. $t\in [0,1]$, it follows
\[
m\{ t \in [0,1] : \gamma'(t) \in \mathfrak{B}  
\text{ and }
\Gamma(t) \neq \gamma(t)\} 
<
\varepsilon.
\]
This proves (2).
\end{proof}

\section{Quotient Maps Preserve Pliability}
\label{s-pliablepush}

\begin{remark}
Suppose $\Phi \colon \mathfrak{g} \to \mathfrak{h}$ is a quotient map with induced Lie group homomorphism $F$. For any horizontal curve $\alpha \colon [0,1] \to \mathbb{G}$,
\begin{equation}
    \label{e-gradedboundlemma}
    \begin{aligned}
    L_\mathbb{H}(F \circ \alpha) 
    =
    \int_0^1 \|(F \circ \alpha)'\|_{\mathfrak{h}_{1}} 
    =
    \int_0^1 \|\Phi \circ (\alpha')\|_{\mathfrak{h}_1} 
    &\leq
    \Vert \Phi|_{\mathfrak{g}_1}\Vert 
    \int_0^1 \|\alpha'\|_{\mathfrak{g}_1} \\
    &=\Vert \Phi|_{\mathfrak{g}_1}\Vert 
    L_\mathbb{G}(\alpha),
\end{aligned}
\end{equation}
where $\Vert \Phi|_{\mathfrak{g}_1}\Vert$ is the operator norm of the linear map $\Phi|_{\mathfrak{g}_1} \colon \mathfrak{g}_1 \to \mathfrak{h}_1$.

As a consequence, we note that $F$ is $\Vert \Phi|_{\mathfrak{g}_1}\Vert$-Lipschitz. Indeed, given $x,y \in \mathbb{G}$ and a horizontal length-minimizing curve $\beta\colon [0,1] \to \mathbb{G}$ from $x$ to $y$, 
the definition of $d_\mathbb{H}$ and \eqref{e-gradedboundlemma} give
\begin{equation}
\label{e-Fis1Lip}
d_\mathbb{H}(F(x),F(y)) \leq L_\mathbb{H}(F \circ \beta) \leq 
\Vert \Phi|_{\mathfrak{g}_1}\Vert
L_\mathbb{G}(\beta) = \Vert \Phi|_{\mathfrak{g}_1}\Vert d_\mathbb{G}(x,y).
\end{equation}
\end{remark}

The following is a bi-Lipschitz alteration of a seemingly well-known result, e.g. see \cite{LPS19}. We give a short proof here.

\begin{lemma}\label{distancelemma}
Suppose $\mathbb{G}, \mathbb{H}$ are Carnot groups and $\Phi\colon \mathfrak{g} \to \mathfrak{h}$ is a quotient map.
Then there is a constant $C=C(\Phi) \geq 1$ such that,
for all $x\in \bbG$ 
and $\tilde{y}\in \mathbb{H}$, there is some $y\in F_\Phi^{-1}(\tilde{y})$ such that
\[
C^{-1} d_{\bbG}(x,y)  \leq d_{\mathbb{H}}(F_\Phi(x), \tilde{y}) \leq C d_{\bbG}(x,y).
\]
\end{lemma}

\begin{proof}
Fix $x\in \bbG$ and $\tilde{y}\in \mathbb{H}$.
Choose a horizontal length-minimizing curve $\gamma\colon [0,1]\to \mathbb{H}$ from $F_\Phi(x)$ to $\tilde{y}$.
Since $F$ is a homomorphism, we may assume that $x = 0_\mathbb{G}$ so that $F_\Phi(x)=0_\mathbb{H}$.
Thus $d_\mathbb{H}(0_\mathbb{H},\tilde{y}) = L_\mathbb{H}(\gamma)$.

Let $\alpha\colon [0,1] \to \mathbb{G}$ be a lift of $\gamma$ by $F_\Phi$ as guaranteed by Lemma~\ref{e-PicardLindeloff} and suppose $C'>0$ is the constant from that lemma. 
Let $y:= \alpha(1)$. Then $F_\Phi(y) = F_\Phi(\alpha(1)) = \gamma(1) = \tilde{y}$,
and
\[
d_{\mathbb{G}}(0_\mathbb{G},y) \leq L_\mathbb{G}(\alpha) \leq C' L_\mathbb{H}(\gamma) = C'd_{\mathbb{H}}(0_\mathbb{H},\tilde{y}).
\]
This proves the first inequality, and the second
follows from \eqref{e-Fis1Lip} since $\tilde{y} = F_\Phi(y)$.
\end{proof}

We now show that quotient mappings preserve pliability.

\begin{proposition}\label{pushforward}
Suppose $\mathbb{G}, \mathbb{H}$ are Carnot groups and $\Phi\colon \mathfrak{g} \to \mathfrak{h}$ is a quotient map.
If $X\in \mathfrak{g}_1$ is pliable, then $\Phi(X)\in \mathfrak{h}_1$ is pliable.
\end{proposition}

\begin{proof}
Using the assumption that $X$ is pliable, 
for every neighborhood $\mathcal{V}$ of the curve $\gamma_X(t) := \exp_{\bbG}(tX)$ in $C_{H}^{1}([0,1],\mathbb{G})$, the set
\[\{\alpha(1):\alpha \in \mathcal{V}, (\alpha, \alpha')(0)=(0_{\bbG},X)\}\]
is a neighborhood of $\gamma_{X}(1)=\exp_{\bbG}(X)$ in $\mathbb{G}$.

Let $Y:=\Phi(X)$ and $\gamma_Y(t) := \exp_{\mathbb{H}}(tY)$ for $t\in [0,1]$. 
By Theorem~\ref{t-fundamental},
\[\exp_{\mathbb{H}}(tY) = \exp_{\mathbb{H}}(\Phi(tX)) = F_\Phi(\exp_{\bbG}(tX)) \qquad \mbox{for every }t\in [0,1].
\]
That is, $\gamma_Y = F_\Phi \circ \gamma_X$.

Fix an arbitrary neighborhood $\mathcal{U}$ of $\gamma_Y$ in $C_{H}^{1}([0,1],\mathbb{H})$. Using the definition of pliability, 
it suffices to show that the set
\[
\{\beta(1):\beta \in \mathcal{U}, (\beta, \beta')(0)=(0_{\mathbb{H}},Y)\}
\]
is a neighborhood of $\gamma_{Y}(1)=\exp_{\mathbb{H}}(Y)$ in $\mathbb{H}$.
To do so, we will show that there is some $s>0$ such that, if $y \in \mathbb{H}$ satisfies $d_{\mathbb{H}}(y,\exp_{\mathbb{H}}(Y)) < s$, then there is a curve $\beta \in \mathcal{U}$ with $(\beta,\beta')(0) = (0_{\mathbb{H}},Y)$ such that $\beta(1)=y$.

Suppose $C \geq \max\{\Vert \Phi|_{\mathfrak{g}_1} \Vert,1\}$ is a constant satisfying the conclusion of Lemma~\ref{distancelemma} for the quotient map $\Phi$.
Using the fact that $\mathcal{U}$ is a neighborhood of $\gamma_Y$, we may choose $r>0$ such that $B(\gamma_Y,r) \subset \mathcal{U}$. Since $X$ is pliable, there is some $0<s<r$ such that
\begin{equation}\label{Xpliable}
B(\exp_{\bbG}(X),s) \subset  \{\alpha(1): \alpha \in B(\gamma_X,\tfrac{r}{C}), \, (\alpha, \alpha')(0)=(0_{\bbG},X)\}.
\end{equation}
Fix $y \in \mathbb{H}$ such that $d(y,\exp_{\mathbb{H}}(Y)) < \frac{s}{C}$.
According to Lemma \ref{distancelemma}, there is some $x \in \mathbb{G}$ with $d(x,\exp_{\mathbb{G}}(X)) \leq Cd(y,\exp_{\mathbb{H}}(Y))<s$ and $F_{\Phi}(x)=y$,
so it follows from \eqref{Xpliable} that there is some $\alpha \in B(\gamma_X,\frac{r}{C})$ with $(\alpha,\alpha')(0)=(0_{\bbG},X)$ and $\alpha(1)=x$.

Consider the curve $\beta = F_{\Phi} \circ \alpha$ in $\mathbb{H}$. We have $\beta(0)=F_{\Phi}(\alpha(0))=F_{\Phi}(0_{\mathbb{G}})=0_{\mathbb{H}}$ since $F_{\Phi}$ is a Lie group homomorphism. Also,
Proposition~\ref{p-chainrule} implies 
\[\beta'(0)=\Phi(\alpha'(0))=\Phi(X)=Y.\]
Hence $(\beta,\beta')(0)=(0_{\mathbb{H}},Y)$. In addition, we note $\beta(1)=F_\Phi(\alpha(1))=F_\Phi(x)=y$. Since $F_\Phi$ is $C$-Lipschitz (see \eqref{e-Fis1Lip}), we have
\[
d_{\mathbb{H}}(\beta(t), \gamma_Y(t)) = d_{\mathbb{H}}(F_\Phi(\alpha(t)),F_\Phi(\gamma_X(t)))\leq Cd_{\mathbb{G}}(\alpha(t),\gamma_X(t)) < r.
\]
Also,
\begin{align*}
\|\beta'(t) - \gamma_Y'(t)\|_{\mathfrak{h}_1} = \|(F_\Phi \circ \alpha)'(t) - (F_\Phi \circ\gamma_X)'(t)\|_{\mathfrak{h}_1} &= \|\Phi(\alpha'(t))-\Phi(\gamma_{X}'(t)) \|_{\mathfrak{h}_1}\\
&\leq \Vert \Phi|_{\mathfrak{g}_1} \|\alpha'(t) - \gamma_X'(t)\|_{\mathfrak{g}_1} < r.
\end{align*}
Hence $\beta \in B(\gamma_Y,r) \subset \mathcal{U}$, and this proves that $Y = \Phi(X)$ is pliable.
\end{proof}

Recall the {\em endpoint map} $\operatorname{End}_\mathbb{G}\colon C^0([0,1],\mathfrak{g}_1) \to \mathbb{G}$ defined as $\operatorname{End}_\mathbb{G}(v)=\gamma(1)$ where $\gamma \in  C^1_H([0,1],\mathbb{G})$ satisfies $\gamma(0)=0_{\bbG}$ and $\gamma'=v$.

\begin{remark}
As noted in \cite{JSS26,JS17}, $X \in \mathfrak{g}_1$ is pliable if and only if the endpoint map $\operatorname{End}_\mathbb{G}$
is open at the constant curve $\delta_X:[0,1] \to \mathfrak{g}_1$ defined as $\delta_X(t) := X$ (we include a short proof of this fact below for completeness).
Thus, to prove the above result, we could show that 
$\operatorname{End}_\mathbb{H}$
is open at the constant curve $\delta_{\Phi(X)}$ 
whenever $\operatorname{End}_\mathbb{G}$ is open at $\delta_{X}$.
Since the map $\Phi$ is open,
it would suffice to prove that 
$\Phi^{-1}(U)$ is a neighborhood of $\delta_{X}$ for any neighborhood $U$ of $\delta_{\Phi(X)}$, i.e., that the projection of $C^0([0,1],\mathfrak{g}_1)$ to $C^0([0,1],\mathfrak{h}_1)$ via composition with $\Phi$ is continuous. This is precisely the content of the above proof.
\end{remark}

Let us now prove the equivalent characterization of pliability just mentioned.

\begin{proposition}
    \label{p-observation}
Suppose $\mathbb{G}$ is a Carnot group.
Then $X \in \mathfrak{g}_1$ is pliable if and only if the map $v \mapsto\operatorname{End}_\mathbb{G}(X+v)$ is open at $0 \in C^0([0,1],\mathfrak{g}_1)$.
\end{proposition}

\begin{proof}
    A vector $X \in \mathfrak{g}_1$ is pliable if,
for any neighborhood $\mathcal{V}$ of $t \mapsto \exp(tX)$ in $C^1_H([0,1],\mathbb{G})$, the set
\[
\{\beta(1) : \beta'(0) = X, \, \beta(0) = 0, \, \beta \in \mathcal{V}\}
\]
is a neighborhood of $\exp(X)$.
This is equivalent to saying that, for any $\varepsilon>0$, the set
\begin{align*}
\{ \operatorname{End}(v) :& \; v \in C^0([0,1],\mathfrak{g}_1), \,
v(0)=X,
\, \Vert v-X\Vert_\infty<\varepsilon\}\\
&= \{\operatorname{End}(X+v): v \in C^0([0,1],\mathfrak{g}_1), \,
v(0)=0,
\, \Vert v\Vert_\infty<\varepsilon\}
\end{align*}
is a neighborhood of $\exp(X)$.
\end{proof}

\section{Applications}
\label{s-applications}

In this section, we use Theorem \ref{t-combine} to investigate pliability and Lusin approximation in several Carnot groups. This is accomplished by construction of suitable quotient mappings onto previously understood spaces.

\subsection{Model-filiform groups}

We first use previous results in the Engel group to infer information about model-filiform groups.

\begin{definition}
\label{d-engelalg}    
The {\em step 3 Engel algebra} $\mathfrak{e}$ is the 4-dimensional Lie algebra with basis $\{E_1,E_2,E_3,E_4\}$ such that the only non-vanishing bracket pairs are
\[
[E_1,E_2] = E_3\quad \text{ and }\quad  [E_1,[E_1,E_2]]=[E_1,E_3]=E_4.
\]
The Engel group $\mathbb{E}$ is the unique Carnot group (up to Lie group isomorphism) whose Lie algebra is $\mathfrak{e}$.
\end{definition}

\begin{definition}
The {\em model filiform algebra} $\mathfrak{fil}_s$ of step $s$ is the $(s+1)$-dimensional Lie algebra  with the stratification
\[
\mathfrak{fil}_s = \operatorname{span}\{Z_1,Z_2\} 
\oplus
\operatorname{span}\{Z_3\}
\oplus
\cdots
\oplus
\operatorname{span}\{Z_{s+1}\}
\]
such that the only non-vanishing bracket pairs are
$[Z_1,Z_k] = Z_{k+1}$ for $k=2,\dots,s$. 
The model Filiform group $\mathbbm{Fil}_{s}$ is the unique Carnot group (up to Lie group isomorphism) whose Lie algebra is $\mathfrak{fil}_s$.
\end{definition}

In particular, $\mathfrak{fil}_1 = \mathbb{R}$, $\mathfrak{fil}_2=\mathfrak{h}$ is the first Heisenberg algebra, and $\mathfrak{fil}_3 = \mathfrak{e}$ is the step 3 Engel algebra. We recall the following fact from \cite{SZ25}.

\begin{proposition}
\label{p-engelnot}
A vector $0 \neq X \in \mathfrak{e}_1$ is pliable if and only if 
$X \notin \operatorname{span}(E_2)$.
\end{proposition}

As a corollary of Proposition \ref{p-engelnot} and Theorem \ref{t-combine}, we obtain the following.

\begin{corollary}
    \label{c-fil}
    For any $s \geq 3$, the vector $Z_2 \in \mathfrak{fil}_s$ is not pliable.
\end{corollary}

\begin{proof}
    Define the Lie algebra homomorphism $\Phi:\mathfrak{fil}_s \to \mathfrak{fil}_3 = \mathfrak{e}$ by extending
    \[
    \Phi(Z_i)=\begin{cases}
    E_i &\text{for } i=1,2,3,4 \\
    0 &\text{for } i \geq 5
    \end{cases}
    \]
    linearly to all of $\mathfrak{fil}_s$. 
    This is indeed a Lie algebra homeomorphism as
    $\Phi([Z_1,Z_2]) = \Phi(Z_3)=E_3=[E_1,E_2]$ and $\Phi([Z_1,Z_3])=[E_1,E_3]$,  and the images of all other bracket pairs vanish. Thus, for any $X = \sum_{i=1}^{s+1} x_i Z_i \in \mathfrak{fil}_s$ and $Y = \sum_{i=1}^{s+1} y_i Z_i \in \mathfrak{fil}_s$, we have
    \begin{align*}
        \Phi([X,Y]) 
        = \sum_{i,j=1}^{s+1} x_iy_j \Phi([Z_i,Z_j])
        = \sum_{i,j=1}^{4} x_iy_j [\Phi(Z_i),\Phi(Z_j)] 
        = [\Phi(X),\Phi(Y)].
    \end{align*}
    Since $\Phi(\operatorname{span}\{Z_1,Z_2\} ) = \operatorname{span}\{E_1,E_2\}$, $\Phi$ is a quotient map. Thus, since Proposition~\ref{p-engelnot} tells us $E_2$ is not pliable, Lemma~\ref{pushforward} implies $Z_2$ cannot be pliable.
\end{proof}

\begin{corollary}
    For any $s \geq 3$, $\mathbbm{Fil}_s$ does not have the $C^1_H$ Lusin property.
\end{corollary}

\begin{proof}
    Recall that $\mathbb{E}$ does not have the $C^1_H$ Lusin property on $\mathfrak{B} = \mathfrak{e}_1$ \cite{Spe16}, so Proposition~\ref{t-surjective} implies that $\mathbbm{Fil}_s$ cannot have the $C^1_H$ Lusin property on $\Phi^{-1}(\mathfrak{B}) \cap (\mathfrak{fil}_s)_1 = (\mathfrak{fil}_s)_1$.
\end{proof}

\subsection{Free Carnot Groups}

In \cite{SZ26} the authors showed that, in the free Carnot group $\mathbb{F}_{2,3}$ of step 3 with 2 generators, every non-zero vector is non-pliable. Further, the weak $C^{1}_{H}$ Lusin property fails, and there exists a Lipschitz curve whose intersection with  {\em any} $C^1_H$ curve has 1-dimensional Hausdorff measure 0. We now apply our earlier work to show this pathological behavior exists in every free Carnot group whose step is at least $3$ and whose number of generators is at least $2$. Recall that the free-nilpotent Lie algebras can be defined as follows \cite[Definition 14.1.1]{BLU07}.

\begin{definition}\label{freeliealgebra}
Let $r\geq 2$ and $s\geq 1$ be integers. We say that $\mathfrak{f}_{r,s}$ is the \emph{free-nilpotent Lie algebra} of \emph{step} $s$ with $r$ \emph{generators} 
if
\begin{enumerate}
\item $\mathfrak{f}_{r,s}$ is a Lie algebra generated by $r$ elements, which we denote by $F^{r,s}_1, \ldots, F^{r,s}_r$,
\item $\mathfrak{f}_{r,s}$ is nilpotent of step $s$ (i.e., nested Lie brackets of length $s+1$ are $0$),
\item for every nilpotent Lie algebra $\mathfrak{g}$ of step $s$ and every map $f\colon \{F^{r,s}_1, \ldots, F^{r,s}_r\}\to \mathfrak{g}$, there is a (unique) Lie algebra homomorphism $\Phi\colon \mathfrak{f}_{r,s} \to \mathfrak{g}$ that extends $f$.
\end{enumerate}
\end{definition}

For each fixed $r\geq 2$ and $s\geq 1$, it can be shown that there exists a Lie algebra $\mathfrak{f}_{r,s}$ satisfying the properties in Definition \ref{freeliealgebra}, and it is unique up to isomorphism of Lie algebras \cite{hall50}.  We next define free Carnot groups \cite[Definition 14.1.3]{BLU07}.

\begin{definition}\label{freecarnotgroup}
Let $r\geq 2$ and $s\geq 1$ be integers. We say that $\mathbb{F}_{r,s}$ is the \emph{free Carnot group} of rank $r$ and step $s$ if its Lie algebra is isomorphic to a free-nilpotent Lie algebra~$\mathfrak{f}_{r, s}$. 
\end{definition}

Note that if two Carnot groups have isomorphic Lie algebras, then the Carnot groups themselves are isomorphic. Hence there is essentially one Carnot group corresponding to each free-nilpotent Lie algebra $\mathfrak{f}_{r, s}$. By saying that two Carnot groups are isomorphic we simply mean that they are isomorphic as Lie groups, with an isomorphism that preserves the stratification.

Let us record the following elementary result.

\begin{lemma}
\label{l-freequot}
For every $r \geq 2$ and $s \geq 3$, 
there is a quotient map $\Phi \colon \mathfrak{f}_{r,s} \to \mathfrak{f}_{2,3}$. 
\end{lemma}

\begin{proof}
Define the map $f:\{F^{r,3}_1,\dots, F^{r,3}_r\} \to \mathfrak{f}_{2,3}$ 
as
\[
f(F^{r,3}_i)=
\begin{cases}
    F^{2,3}_i &\text{if } i\in \{1,2\}\\
    0 &\text{if } i  \in \{3,\dots,r\}
\end{cases}
\]
By definition, there is a unique Lie algebra homomorphism $\Phi_1:\mathfrak{f}_{r,3} \to \mathfrak{f}_{2,3}$ extending $f$.

By identifying $F_i^{r,s} \leftrightarrow F_i^{r,3}$ for $i \in \{1,\dots,r\}$, 
we may identify $(\mathfrak{f}_{r,s})_k = (\mathfrak{f}_{r,3})_k$ for $k=1,2,3$ 
since both $\mathfrak{f}_{r,s}$ and $\mathfrak{f}_{r,3}$ are Lie algebras generated by these basis elements.
Thus, defining $\Phi_2:\mathfrak{f}_{r,s} \to \mathfrak{f}_{r,3}$ to be the identity on $(\mathfrak{f}_{r,s})_1 \oplus (\mathfrak{f}_{r,s})_2 \oplus (\mathfrak{f}_{r,s})_3$ and setting $\Phi_2\left((\mathfrak{f}_{r,s})_k \right) = \{0\}$ for any $k \geq 4$ ensures that $\Phi_2$ is a quotient map. 
In particular, $\Phi_1(\Phi_2((\mathfrak{f}_{r,s})_1)) = \Phi_1((\mathfrak{f}_{r,3})_1) = (\mathfrak{f}_{2,3})_1$.
    Defining $\Phi = \Phi_1 \circ \Phi_2$ completes the proof.
\end{proof}

The following shows that the Lusin property fails in a strong sense in free Carnot groups.

\begin{corollary}    
    Fix integers $r \geq 2$ and $s \geq 3$.
    Then there is a Lipschitz curve $\gamma \colon [0,1] \to \mathbb{F}_{r,s}$ such that
    \[
    m\{t \in [0,1] : \gamma(t) = \Gamma(t) \} = 0
    \]
    for any $C^1_H$ curve $\Gamma \colon [0,1] \to \mathbb{F}_{r,s}$.
\end{corollary}

\begin{proof}
Suppose the statement does not hold. Then, by definition, $\mathbb{F}_{r,s}$ has the weak $C^1_H$ Lusin property. By Lemma~\ref{l-freequot} there is a quotient map from $\mathfrak{f}_{r,s}$ to $\mathfrak{f}_{2,3}$, so Proposition~\ref{t-surjective}  implies that $\mathfrak{f}_{2,3}$ has the weak $C^1_H$ Lusin property. However, this contradicts Theorem 1.2 from \cite{SZ26}.
\end{proof}

Combining this result with Theorem~5.1 from \cite{SZ26}, we have the following infinite collection of pathological examples.
Recall that, in a Carnot group $\mathbb{G}$, a set $E \subset \mathbb{G}$ is purely unrectifiable by $C^1_H$ curves if 
\[
\mathcal{H}^1(E \cap \Gamma([0,1])) = 0
\]
for every $C^1_H$ curve $\Gamma:[0,1] \to \mathbb{G}$.
Here, $\mathcal{H}^1$ is standard 1-dimensional Hausdorff measure with respect to the CC metric.

\begin{corollary}
        Fix integers $r \geq 2$ and $s \geq 3$.
    Then there is a Lipschitz curve $\gamma\colon [0,1] \to \mathbb{F}_{r,s}$ such that
    $\gamma([0,1])$ is purely unrectifiable by $C^1_H$ curves.
\end{corollary}

Using similar techniques, we show that every non-zero vector is non-pliable.

\begin{lemma}\label{maptoZ2}
Fix $0 \neq X \in (\mathfrak{f}_{r,s})_1$.
    There is a quotient map $\Phi \colon \mathfrak{f}_{r,s} \to \mathfrak{fil}_{s}$ with $\Phi(X) = Z_2$.
\end{lemma}

\begin{proof}
Let $F_1,\dots,F_r$ denote the generators of $\mathfrak{f}_{r,s}$, and write $X = \sum_{k=1}^r x_kF_k \neq 0$ for some $x_{1}\ldots, x_{r}\in \mathbb{R}$.
To prove the lemma, we will define a map $\phi \colon \{F_1,\dots,F_r\} \to \mathfrak{fil}_s$ in such a way that, when $\phi$ is extended to a Lie algebra homomorphism $\Phi$ as in Definition~\ref{freeliealgebra}, we will have $\Phi(X)=Z_2$ and $\Phi$ will be a quotient map. 

Without loss of generality, assume $x_1 \neq 0$.
    Define 
    $\phi \colon \{F_1,\dots,F_r\} \to \mathfrak{fil}_s$ so that
    $\phi(F_1)=\frac{1}{x_1}Z_2 - \frac{x_2}{x_1}(Z_1+Z_2)$
    and
    $\phi(F_2)=Z_1+Z_2$, and set $\phi(F_k)=0$ for $k \geq 3$.
    By Definition~\ref{freeliealgebra}, 
    there is a Lie algebra homomorphism $\Phi:\mathfrak{f}_{r,s} \to \mathfrak{fil}_s$ extending $\phi$. In particular, 
    \[
    \Phi(-x_1F_1 + (1-x_2)F_2) = - Z_2 + x_2(Z_1+Z_2) + (1-x_2)(Z_1 + Z_2) 
    = Z_1
    \]
    and
    \[
    \Phi(X) =  x_1\Phi(F_1)+ x_2\Phi(F_2) = Z_2 - x_2(Z_1+Z_2) +x_2(Z_1 + Z_2) 
    =Z_2.
    \]
    Hence $\Phi((\mathfrak{f}_{r,s})_1) = (\mathfrak{fil}_s)_1$, so $\Phi$ is a quotient map. Since also $\Phi(X)=Z_{2}$, the proof is complete.
\end{proof}

\begin{corollary}
    The only pliable vector in $(\mathfrak{f}_{r,s})_1$ is $0$.
\end{corollary}

\begin{proof}
Suppose $0 \neq X \in (\mathfrak{f}_{r,s})_1$ is pliable. Using Lemma \ref{maptoZ2}, choose a quotient map $\Phi \colon \mathfrak{f}_{r,s} \to \mathfrak{fil}_{s}$ with $\Phi(X) = Z_2$. Applying Lemma~\ref{pushforward} would imply that $Z_{2}$ is pliable, but this contradicts Corollary~\ref{c-fil}. Hence $X$ cannot be pliable.
\end{proof}

\subsection{Nilpotent matrices}

For positive integer $r$, let $\mathbb{U}_{r}$ be the Lie group of ${r} \times {r}$ upper triangular {\em unipotent} matrices with real valued entries, i.e., $r \times r$ upper triangular matrices with 1's on the diagonal. This is a Carnot group of step $(r-1)$ \cite[Example 9.3.4]{MR4995404}. It's Lie algebra $\mathfrak{n}_r$ is the vector space of strictly upper triangular {\em $r$-nilpotent} real valued matrices, i.e., $r \times r$ upper triangular matrices with 0's on the diagonal \cite[Example 9.1.9]{MR4995404}.

Denote by $E_{i,j}$ the matrix with a 1 in the $i$th row and $j$th column and 0 everywhere else.
Then $\{ E_{i,j}: 1 \leq i < j \leq r\}$ is a basis of $\mathfrak{n}_r$.
Moreover, 
$\mathfrak{n}_r$ is a step $(r-1)$ nilpotent Lie algebra with the stratification 
\[
(\mathfrak{n}_r)_k = \operatorname{span}\{E_{i,i+k}:i\in\{1,\dots,r-k\}\} \quad
\text{for $k \in \{1,\dots,r-1\}$}
\]
since we have $[E_{i,i+1},E_{i+1,i+k+1}]=E_{i,i+k+1}$ for all $i \in \{1,\dots,r-k-1\}$. 

\begin{corollary}
    \label{c-nipotentexample}
    For each positive integer $r$, $\mathbb{U}_r$ is pliable. 
\end{corollary}

To prove this, we will construct a quotient map from the Lie algebra of a pliable Carnot group to the Lie algebra of $\mathbb{U}_r$.
Given an integer $s \geq 3$, we follow the example of \cite[Proposition~6.6]{JS17} to define a step $s$, stratified, nilpotent Lie algebra $\mathfrak{G}_s$ whose associated Carnot group is pliable. Suppose $F_1,\dots,F_s$ are the generators of $\mathfrak{f}_{s,s}$, and define $I_i$ to be the ideal in $\mathfrak{f}_{s,s}$ generated by $F_i$ for each $i \in \{1,\dots,s\}$. 
That is, $I_i$ is the smallest linear subspace of $\mathfrak{f}_{s,s}$ containing $F_i$ and satisfying $[X,Y] \in I_i$ for all $X \in I_i$ and $Y \in \mathfrak{f}_{s,s}$.
Then $J = \sum_{i=1}^s [I_i,I_i]$ is also an ideal, so the Lie algebra 
$\mathfrak{G}_s=\mathfrak{f}_{s,s}/J$
is also stratified and nilpotent (see, for example, \cite[Exercise 9.5.14]{MR4995404}).
According to the proof of Proposition~6.6 in \cite{JS17}, the Carnot group associated with $\mathfrak{G}_s$ is pliable.

\begin{example}
One may check that $\mathfrak{G}_3$ satisfies 
\[
\mathfrak{G}_3=
(\mathfrak{G}_3)_1
\oplus
(\mathfrak{G}_3)_2
\oplus
(\mathfrak{G}_3)_3
\]
with
\[
\mathfrak{G}^{(3)}_1 = \operatorname{span}\{X_1,X_2,X_3\},
\qquad
\mathfrak{G}^{(3)}_2 = \operatorname{span}\{X_{12},X_{13},X_{23}\},
\qquad
\mathfrak{G}^{(3)}_3 = \operatorname{span}\{X_{123},X_{213}\}.
\]
Here, the non-vanishing bracket relations are
\[
[X_1,X_2] = X_{12},
\qquad
[X_1,X_3] = X_{13},
\qquad
[X_2,X_3] = X_{23}
\]
\[
[X_1,X_{23}] = X_{123},
\qquad
[X_2,X_{13}] = X_{213}
\]
up to skew-symmetry,
and $[X_3,X_{12}] = X_{213}-X_{123}$ (where this last bracket relation is required by the Jacobi identity).
That is, $\mathfrak{G}_3$ is an 8-dimensional, stratified, nilpotent Lie algebra of step 3.
\end{example}

\begin{proof}[Proof of Corollary \ref{c-nipotentexample}]
Fix a positive integer $s$, and set $r = s+1$.
If $s=1$ or 2, then $\mathbb{U}_{r}$ is pliable since
it has step at most 2 (see \cite[Theorem~6.5]{JS17}).
Suppose $s \geq 3$,
and define $\mathfrak{G}_s$ as above.
For convenience, we will use the identification
$(\mathfrak{G}_s)_1
=
(\mathfrak{f}_{s,s})_1 = \operatorname{span}\{F_1,\dots,F_s\}$.

To prove that every horizontal vector in $\mathfrak{n}_{s+1}$ is pliable, we will construct a quotient map $\Phi\colon \mathfrak{G}_s \to \mathfrak{n}_{s+1}$.
First, define $\tilde{\Phi}(F_i)=E_{i,i+1}$ for $i = 1,\dots,s$. By definition, this extends to a Lie algebra homomorphism $\tilde{\Phi}\colon \mathfrak{f}_{s,s} \to \mathfrak{n}_{s+1}$.
We will show that this induces a Lie algebra homomorphism $\Phi\colon \mathfrak{G}_s \to \mathfrak{n}_{s+1}$ by proving the following claim.
Recall that we defined $J = \sum_{i=1}^s [I_i,I_i]$ where $I_i$ is the ideal in $\mathfrak{f}_{s,s}$ generated by $F_i$.

\medskip

\noindent
\textbf{Claim:} $J \subset \operatorname{ker}\tilde{\Phi}$.

\begin{proof}[Proof of Claim]
It suffices to prove that 
$
\tilde{\Phi}([I_i,I_i]) = [\tilde{\Phi}(I_i),\tilde{\Phi}(I_i)] = \{0\},
$
and this is simply a linear algebraic argument.
Suppose $M$ is an $r\times r$ matrix of the form 
\begin{equation}
\label{e-blockupperright}
M=\left( 
\begin{array}{cc}
   0  &  A\\
   0  & 0
\end{array}\right)
\end{equation}
for some $i \times (r-i)$ matrix $A$.
Then, if $N$ is any $r \times r$ upper triangular
matrix, it follows that $MN$ and $NM$ both have the form \eqref{e-blockupperright}.
Hence the nested bracket of $M$ with any $r \times r$ upper triangular matrices must have the form \eqref{e-blockupperright}. 

Fix $X \in \tilde{\Phi}(I_i)$. Then $X$ is equal to a linear combination of nested brackets of $\tilde{\Phi}(F_i)=E_{i,i+1}$ with (upper triangular) elements of $\mathfrak{n}_r$. Since $E_{i,i+1}$ is of the form \eqref{e-blockupperright}, $X$ must have the form \eqref{e-blockupperright} as well.
Note, however, that $M_1M_2=0$ for any $r \times r$ matrices of the form \eqref{e-blockupperright}, so it follows that $[\tilde{\Phi}(I_i),\tilde{\Phi}(I_i)] = \{0\}$ as desired.
\end{proof}
That is, there is a Lie algebra homomorphism $\Phi \colon \mathfrak{G}_s \to \mathfrak{n}_{s+1}$.
It remains to prove that $\Phi$ is a quotient map, and this follows from the definition of $\tilde{\Phi}$ since
\[
\Phi((\mathfrak{G}_s)_1)
=
\tilde{\Phi}((\mathfrak{f}_{s,s})_1)
=
\operatorname{span}\{E_{12},E_{23},\dots,E_{s,s+1}\}=(\mathfrak{n}_{s+1})_1.
\]

\end{proof}

\printbibliography

@article {Zim18,
    AUTHOR = {Zimmerman, Scott},
     TITLE = {The {W}hitney extension theorem for {$C^1$}, horizontal curves
              in the {H}eisenberg group},
   JOURNAL = {J. Geom. Anal.},
  FJOURNAL = {Journal of Geometric Analysis},
    VOLUME = {28},
      YEAR = {2018},
    NUMBER = {1},
     PAGES = {61--83},
      ISSN = {1050-6926,1559-002X},
   MRCLASS = {53C17},
  MRNUMBER = {3745849},
MRREVIEWER = {Matthew\ Romney},
       DOI = {10.1007/s12220-017-9807-2},
       URL = {https://doi-org.proxy.lib.ohio-state.edu/10.1007/s12220-017-9807-2},
}

@article {JS17,
    AUTHOR = {Juillet, Nicolas and Sigalotti, Mario},
     TITLE = {Pliability, or the {W}hitney extension theorem for curves in
              {C}arnot groups},
   JOURNAL = {Anal. PDE},
  FJOURNAL = {Analysis \& PDE},
    VOLUME = {10},
      YEAR = {2017},
    NUMBER = {7},
     PAGES = {1637--1661},
      ISSN = {2157-5045,1948-206X},
   MRCLASS = {22E25 (53C17 54C20 58C25 93B27)},
  MRNUMBER = {3683924},
MRREVIEWER = {Andrea\ Pinamonti},
       DOI = {10.2140/apde.2017.10.1637},
       URL = {https://doi-org.proxy.lib.ohio-state.edu/10.2140/apde.2017.10.1637},
}

@article {SZ25,
    AUTHOR = {Speight, Gareth and Zimmerman, Scott},
     TITLE = {Directional pliability, {W}hitney extension, and {L}usin
              approximation for curves in {C}arnot groups},
   JOURNAL = {Ann. Fenn. Math.},
  FJOURNAL = {Annales Fennici Mathematici},
    VOLUME = {50},
      YEAR = {2025},
    NUMBER = {2},
     PAGES = {665--684},
      ISSN = {2737-0690,2737-114X},
   MRCLASS = {53C17 (58C25)},
  MRNUMBER = {4977317},
       DOI = {10.54330/afm.176524},
       URL = {https://doi-org.proxy.lib.ohio-state.edu/10.54330/afm.176524},
}

@article {Spe16,
    AUTHOR = {Speight, Gareth},
     TITLE = {Lusin approximation and horizontal curves in {C}arnot groups},
   JOURNAL = {Rev. Mat. Iberoam.},
  FJOURNAL = {Revista Matem\'atica Iberoamericana},
    VOLUME = {32},
      YEAR = {2016},
    NUMBER = {4},
     PAGES = {1423--1444},
      ISSN = {0213-2230,2235-0616},
   MRCLASS = {53C17 (26A46)},
  MRNUMBER = {3593531},
MRREVIEWER = {Andrea\ Pinamonti},
       DOI = {10.4171/RMI/924},
       URL = {https://doi-org.proxy.lib.ohio-state.edu/10.4171/RMI/924},
}

@book {BLU07,
    AUTHOR = {Bonfiglioli, A. and Lanconelli, E. and Uguzzoni, F.},
     TITLE = {Stratified {L}ie groups and potential theory for their
              sub-{L}aplacians},
    SERIES = {Springer Monographs in Mathematics},
 PUBLISHER = {Springer, Berlin},
      YEAR = {2007},
     PAGES = {xxvi+800},
%      ISBN = {978-3-540-71896-3; 3-540-71896-6},
   MRCLASS = {22E30 (31C45 35-02 35H10 43A80)},
  MRNUMBER = {2363343},
MRREVIEWER = {Maria\ Stella\ Fanciullo},
}

@article {LPS19,
    AUTHOR = {Le Donne, Enrico and Pinamonti, Andrea and Speight, Gareth},
     TITLE = {Universal differentiability sets and maximal directional
              derivatives in {C}arnot groups},
   JOURNAL = {J. Math. Pures Appl. (9)},
  FJOURNAL = {Journal de Math\'ematiques Pures et Appliqu\'ees. Neuvi\`eme
              S\'erie},
    VOLUME = {121},
      YEAR = {2019},
     PAGES = {83--112},
      ISSN = {0021-7824,1776-3371},
   MRCLASS = {53C17 (22E25)},
  MRNUMBER = {3906166},
MRREVIEWER = {Scott\ Robert\ Zimmerman},
       DOI = {10.1016/j.matpur.2017.11.006},
       URL = {https://doi-org.proxy.lib.ohio-state.edu/10.1016/j.matpur.2017.11.006},
}

@article {Pan89,
    AUTHOR = {Pansu, Pierre},
     TITLE = {M\'etriques de {C}arnot-{C}arath\'eodory et quasiisom\'etries
              des espaces sym\'etriques de rang un},
   JOURNAL = {Ann. of Math. (2)},
  FJOURNAL = {Annals of Mathematics. Second Series},
    VOLUME = {129},
      YEAR = {1989},
    NUMBER = {1},
     PAGES = {1--60},
      ISSN = {0003-486X,1939-8980},
   MRCLASS = {53C20 (22E40)},
  MRNUMBER = {979599},
MRREVIEWER = {Gudlaugur\ Thorbergsson},
       DOI = {10.2307/1971484},
       URL = {https://doi-org.proxy.lib.ohio-state.edu/10.2307/1971484},
}

@article {Whi34,
    AUTHOR = {Whitney, Hassler},
     TITLE = {Analytic extensions of differentiable functions defined in
              closed sets},
   JOURNAL = {Trans. Amer. Math. Soc.},
  FJOURNAL = {Transactions of the American Mathematical Society},
    VOLUME = {36},
      YEAR = {1934},
    NUMBER = {1},
     PAGES = {63--89},
      ISSN = {0002-9947,1088-6850},
   MRCLASS = {26B05},
  MRNUMBER = {1501735},
       DOI = {10.2307/1989708},
       URL = {https://doi-org.proxy.lib.ohio-state.edu/10.2307/1989708},
}

@unpublished{SZ26,
  author    = {Gareth Speight and Scott Zimmerman},
  title     = {A Lipschitz curve in a Carnot group that is purely unrectifiable by smooth horizontal curves},
  year      = {2026},
  note      = {Submitted for publication}
}

@article {JSS26,
    AUTHOR = {Jean, Fr\'ed\'eric and Sigalotti, Mario and Socionovo,
              Alessandro},
     TITLE = {A note on pliability and the openness of the multiexponential
              map in {C}arnot groups},
   JOURNAL = {SIGMA Symmetry Integrability Geom. Methods Appl.},
  FJOURNAL = {SIGMA. Symmetry, Integrability and Geometry. Methods and
              Applications},
    VOLUME = {22},
      YEAR = {2026},
     PAGES = {Paper No. 078, 8},
      ISSN = {1815-0659},
   MRCLASS = {53C17 (22E25 93B05)},
  MRNUMBER = {5118634},
       DOI = {10.3842/SIGMA.2026.078},
       URL = {https://doi-org.proxy.lib.ohio-state.edu/10.3842/SIGMA.2026.078},
}

@article {LDS16,
    AUTHOR = {Le Donne, Enrico and Speight, Gareth},
     TITLE = {Lusin approximation for horizontal curves in step 2 {C}arnot
              groups},
   JOURNAL = {Calc. Var. Partial Differential Equations},
  FJOURNAL = {Calculus of Variations and Partial Differential Equations},
    VOLUME = {55},
      YEAR = {2016},
    NUMBER = {5},
     PAGES = {Art. 111, 22},
      ISSN = {0944-2669,1432-0835},
   MRCLASS = {28C15 (43A80 49Q15 53C17)},
  MRNUMBER = {3549917},
MRREVIEWER = {Davide\ Vittone},
       DOI = {10.1007/s00526-016-1054-z},
       URL = {https://doi-org.proxy.lib.ohio-state.edu/10.1007/s00526-016-1054-z},
}

@book {rifford,
    AUTHOR = {Rifford, Ludovic},
     TITLE = {Sub-{R}iemannian geometry and optimal transport},
    SERIES = {SpringerBriefs in Mathematics},
 PUBLISHER = {Springer, Cham},
      YEAR = {2014},
     PAGES = {viii+140},
      ISBN = {},
   MRCLASS = {49Q20 (53C17)},
  MRNUMBER = {3308395},
MRREVIEWER = {Yanir\ A.\ Rubinstein},
       DOI = {10.1007/978-3-319-04804-8},
       URL = {https://doi-org.proxy.lib.ohio-state.edu/10.1007/978-3-319-04804-8},
}

@book {Lee13,
    AUTHOR = {Lee, John M.},
     TITLE = {Introduction to smooth manifolds},
    SERIES = {Graduate Texts in Mathematics},
    VOLUME = {218},
   EDITION = {Second},
 PUBLISHER = {Springer, New York},
      YEAR = {2013},
     PAGES = {xvi+708},
      ISBN = {978-1-4419-9981-8},
   MRCLASS = {58-01 (53-01 57-01)},
  MRNUMBER = {2954043},
}

@article {hall50,
    AUTHOR = {Hall, Jr., Marshall},
     TITLE = {A basis for free {L}ie rings and higher commutators in free
              groups},
   JOURNAL = {Proc. Amer. Math. Soc.},
  FJOURNAL = {Proceedings of the American Mathematical Society},
    VOLUME = {1},
      YEAR = {1950},
     PAGES = {575--581},
      ISSN = {0002-9939,1088-6826},
   MRCLASS = {09.1X},
  MRNUMBER = {38336},
MRREVIEWER = {K.\ A.\ Hirsch},
       DOI = {10.2307/2032282},
       URL = {https://doi-org.proxy.lib.ohio-state.edu/10.2307/2032282},
}

@book {MR4995404,
    AUTHOR = {Le Donne, Enrico},
     TITLE = {Metric {L}ie groups---{C}arnot-{C}arath\'eodory spaces from
              the homogeneous viewpoint},
    SERIES = {Graduate Texts in Mathematics},
    VOLUME = {306},
 PUBLISHER = {Springer, Cham},
      YEAR = {2025},
     PAGES = {xvi+480},
   MRCLASS = {22E25 (22-01 22E30 22F30 53C17)},
  MRNUMBER = {4995404},
MRREVIEWER = {Vladimir\ V.\ Gorbatsevich},
       DOI = {10.1007/978-3-031-98832-5},
       URL = {https://doi-org.uc.idm.oclc.org/10.1007/978-3-031-98832-5},
}

\end{document}